\documentstyle[11pt,leqno]{article}
\newtheorem{theorem}{{\sc Theorem}}
\newcommand{\bt}{\begin{theorem}}
\newcommand{\et}{\end{theorem}}
\newcommand{\newsection}[1]{\setcounter{equation}{0} \setcounter{theorem}{0}
\section{#1}}

\newcommand{\NI}{\noindent}
\newcommand{\bea}{\begin{eqnarray}}
\newcommand{\eea}{\end{eqnarray}}

\def \spec#1 {\mathop{#1}}

\def \b #1 {\bf #1}

\newcommand {\CC}{\centerline}

\newcommand{\clf}{{\cal F}}

\newcommand{\ity}{\infty}
\newcommand{\raro}{\rightarrow}

\newcommand{\vsp}{\vskip 1em}

\newcommand{\ve}{\varepsilon}

\newcommand{\be}{\begin{equation}}
\newcommand{\ee}{\end{equation}}
\newcommand{\ben}{\begin{eqnarray*}}
\newcommand{\een}{\end{eqnarray*}}

\begin{document}

\sloppy
\CC {\Large {\bf Nonparametric Estimation of Linear Multiplier in}}
\CC{\Large {\bf SDEs driven by General Gaussian Processes}} 
\vsp
\CC {\bf B.L.S. Prakasa Rao   }
\CC{\bf CR Rao Advanced Institute of Mathematics, Statistics }
\CC{\bf and Computer Science, Hyderabad, India}
\vsp
\NI{\bf Abstract:} We investigate the asymptotic properties of a kernel-type nonparametric estimator of the linear multiplier in models governed by a stochastic differential equation driven by a general Gaussian process.
\vsp
\newsection{Introduction}

Diffusion processes and diffusion type processes satisfying stochastic
differential equations driven by Wiener processes are used for stochastic
modeling in a wide variety of sciences such as population genetics,
economic processes, signal processing as well as for modeling sunspot
activity and more recently in mathematical finance. Statistical inference
for diffusion type processes satisfying stochastic  differential equations
driven by Wiener processes have been studied earlier and a comprehensive
survey of various methods is given in Prakasa Rao (1999). There has been
a recent interest to study similar problems for stochastic processes
driven by a fractional Brownian motion to model processes involving long
range dependence (cf. Prakasa Rao (2010)).

Our aim in this paper is to investigate the asymptotic properties of a non-parametric kernel-type estimator for estimating the linear multiplier in a stochastic differential equation driven by a Gaussian process.  Nonparametric estimation of a linear multiplier for processes driven by $\alpha$-stable noise is investigated in Prakasa Rao (2021). Asymptotic properties of  the minimum $L_1$-norm estimator of the drift parameter of a Ornstein-Uhlenbeck  process driven by general Gaussian processes are investigated in Prakasa Rao (2022). El Machkouri et al. (2015), Chen and Zhou (2020) and Lu (2022) study parameter estimation for an Ornstein-Uhlenbeck process  driven by a general Gaussian process.

\newsection{Preliminaries}
Let $(\Omega, \clf, (\clf_t), P) $ be a stochastic basis satisfying the
usual conditions and the processes discussed in the following are
$(\clf_t)$-adapted. Further the natural filtration of a process is
understood as the $P$-completion of the filtration generated by this
process. We consider a centered Gaussian process $G\equiv\{G_t,0\leq t \leq T\}.$
\vsp
Let us consider a stochastic process $\{X_t, t \in [0,1]\}$ defined by the 
stochastic differential  equation 
\be
dX_t= \theta(t) X_t dt  + \ve \; dG_t,X_0=x_0 0 \leq t \leq 1,
\ee
where the function $\theta (t)$ is an unknown. We assume that the Gaussian process $G(.)$ satisfies the condition $E[\sup_{0\leq s \leq T}|G(s)|]$ is finite and it has Holder continous paths of positive order. We assume that integration with respect to the Gaussian process $G$ is defined as a Young integral (cf. Nourdin (2012)). This class of Gaussian processes $G$ includes fractional Brownian motion, sub-fractional Brownian motion and bifractional Brownian motion (cf. Mishura and Zili (2018)) under some conditions.
\vsp
We now consider the problem of estimation of the function $\theta(t) $ based
on the observation of process $ X=\{X_t , 0 \leq t \leq 1\}$  and study its
asymptotic properties as $\ve \raro 0.$ 
\vsp
Let $x= \{x_t, 0 \leq t \leq 1\}$ be the solution of the ordinary differential equation
\be
\frac{dx_t}{dt}=\theta(t) x_t , x_0 , 0 \leq t \leq 1.
\ee
Observe that
$$x_t= x_0 \exp(\int_0^t\theta(s)ds).$$
\vsp
We assume that (A1) the function $\theta(t)$ is bounded over the interval $[0,T]$ by a constant $L$. 
\vsp
\NI{\bf Lemma 2.1: } {\it Let  $ X_t$ and $x_t$ be  the solutions of the equation (2.1) and
(2.2) respectively. Then, with probability one,}
\be
(a)|X_t-x_t|\leq  e^{Lt} \epsilon \sup_{0\leq s \leq t}|G_s|
\ee
and
\be
(b)\sup_{0 \leq t \leq T} E|X_t-x_t| \leq e^{LT} \epsilon  E[\sup_{0\leq t \leq T}|G(t)|].
\ee
\vsp
\NI{\bf Proof of (a) :} Let $u_t=|X_t-x_t| $. Then
\bea
u_t & \leq & \int^t_0 \left| \theta(v)(X_v-x_v) \right| dv + \epsilon \;|G_t|\\\nonumber
& \leq & L \int^t_0 u_v dv + \epsilon \;\sup_{0\leq s \leq t}|G_s|.
\eea
Applying the Gronwall's lemma (cf. Lemma 1.12, Kutoyants (1994), p.26),  it follows that
\be
u_t \leq  \epsilon e^{LT}\sup_{0\leq s \leq t}|G_s| . 
\ee
\vsp
\NI{\bf Proof of (b) :} From (2.3), we have ,
\be
E|X_t-x_t| \leq e^{Lt} \epsilon \; E (|\sup_{0\leq s \leq t}|G_s|). 
\ee
Hence
\be
\sup_{0 \leq t \leq T} E |X_t-x_t|  \leq \epsilon e^{LT}E[\sup_{0\leq t \leq T}|G(t)|]. \\
\ee
\vsp
\newsection{Estimation of the Drift function}

Let $\Theta_0(L)$ denote the class of all functions $\theta(.) $ with the same bound $L$. Let $\Theta_k(L) $ denote the class of all functions $\theta(.) $ which are uniformly bounded by the same constant $L$ and which are $k$-times differentiable satisfying the condition
$$|\theta^{(k)}(x)-\theta^{(k)}(y)|\leq L^\prime |x-y|, x,y \in R$$
for some constant $L^\prime >0.$ Here $g^{(k)}(x)$ denotes the $k$-th derivative of $g(.)$ at $x$ for $k \geq 0.$ If $k=0,$ we interpret the function $g^{(0)}(.)$ as the function  $g(.).$\\
Let $K(u)$ be a  bounded function  with finite support $[A,B]$ with $A<0<B$ satisfying the condition 

\NI{$(A_2)$}$ K(u) =0\;\; \mbox{for}\;\; u <A \;\;\mbox{and}\;\; u>B;\;\;\mbox{and} \int^B_A K(u) du =1.$

Boundedness of the function $K(.)$ with finite support $[A,B]$ implies that

$$\int_{-\ity}^\ity |K(u)|^2 du <\ity; \int_{-\ity}^\ity |u^jK(u)| du <\ity, j \geq 0.$$

We define a kernel type estimator $\widehat{\theta}_t$ of the function $\theta(t)$ by the relation 
\be
\widehat{\theta}_t X_t = \frac{1}{\varphi_\epsilon}\int^T_0 K\left(\frac{\tau-t}{\varphi_\epsilon} \right) d X_\tau
\ee
where the normalizing function  $ \varphi_\epsilon \rightarrow 0 $ as  $ \epsilon \rightarrow 0. $ Let $E_\theta(.)$ denote the expectation when the function $\theta(.)$ is the linear multiplier.
\vsp
\NI{\bf Theorem 3.1:}  {\it Suppose that the linear multiplier $\theta(.) \in \Theta_0(L)$ and  the function
\ $ \varphi_\epsilon \rightarrow 0$  as $\epsilon \rightarrow  0$. Further suppose that the conditions $(A_1), (A_2)$   hold. Then, for any $ 0 < c \leq d < T ,$ the estimator $\widehat{\theta}_t$ satisfies the property}
\be
\lim_{\epsilon \rightarrow 0} \sup_{\theta(.) \in \Theta_0(L)} \sup_{c\leq t \leq d } E_\theta ( |\widehat{\theta}_t X_t - \theta(t)x_t|)= 0.
\ee
\vsp
In addition to the conditions $(A_1),(A_2),$ suppose the following condition holds.\\
\NI{$(A_3)$}$ \int^\infty_{-\infty} u^j K(u)   du = 0 \;\;\mbox{for}\;\; j=1,2,...k.$
\vsp
\NI{\bf Theorem 3.2:} {\it Suppose that the function $ \theta(.) \in \Theta_{k+1}(L)$ and  $
\varphi_\epsilon = \epsilon^{\frac{1}{k+1}}.$ Suppose the conditions $(A_1)-(A_3)$ hold. Then}
\be
\limsup_{\epsilon \rightarrow 0} \sup_{\theta(.) \in \Theta_{k+1}(L)}\sup_{c \leq t \leq d} E_\theta (| \widehat{\theta}_t X_t- \theta(t)x_t|)\epsilon^{-1}  < \infty.
\ee
\vsp
\NI{\bf Theorem 3.3:} {\it Suppose that the function $\theta(.) \in \Theta_{k+1}(L)$ and $ \varphi_\epsilon= \epsilon^{1/(k+1)}$. Further suppose that the conditions $(A_1)-(A_3)$ hold. Then the asymptotic distribution of
$$ \varphi_\epsilon^{-(k+1)}(\widehat{\theta}_t X_t- \theta(t)x_t), $$
as $\epsilon \raro 0$ is the distribution of a Gaussian  random variable with mean
$$m= \frac{J^{(k+1)}(x_t)}{(k+1)!}\int_{-\ity}^{\ity}K(u)u^{k+1}du$$
and variance $R(t,t)$ as $\epsilon \raro 0$ where $J(t)=\theta(t)x(t).$}
\vsp
\newsection{Proofs of Theorems 3.1-3.3:}

\NI{\bf Proof of Theorem 3.1 :} From the equation (3.1), we have
\bea
\;\;\;\\\nonumber
E_\theta[|\widehat{\theta}_tX_t -\theta(t)x_t|]  &=& E_\theta [ |\frac{1}{\varphi_\epsilon}  \int^T_0 K \left(\frac{\tau-t}{\varphi_\epsilon} \right) \left(\theta(\tau)X_\tau -\theta(\tau) x_\tau \right)  d \tau \\ \nonumber
& &+ \frac{1}{\varphi_\epsilon}
 \int^T_0 K\left(\frac{\tau-t}{\varphi_\epsilon}\right) \theta(\tau) x_\tau d \tau- \theta(t) x_t
 + \frac{\epsilon}{\varphi_\epsilon} \int^T_0 K \left(\frac{\tau-t}{\varphi_\epsilon} \right)
 dG_\tau|]\\ \nonumber
 & \leq  &  E_\theta \left[ |\frac{1}{\varphi_\epsilon}  \int^T_0 K \left(\frac{\tau-t}{\varphi_\epsilon} \right) (\theta(\tau)X_\tau -\theta(\tau)x_\tau) d\tau |\right] \\ \nonumber
 & & + E_\theta \left[|\frac{1}{\varphi_\epsilon} \int^T_0 K \left(\frac{\tau-t}{\varphi_\epsilon} \right)\theta(\tau) x_\tau d\tau -\theta(t) x_t |\right] \\\nonumber
 & & +  \frac{\epsilon}{\varphi_{\epsilon}} E_\theta \left[ |\int^T_0 K \left(  \frac{\tau-t}{\varphi_\epsilon}\right) d G_\tau|\right]\\ \nonumber
 &= & I_1+I_2+I_3 \;\;\mbox{(say).}\;\;\\ \nonumber
\eea
Apply the change of variables $u= (\tau-t)\varphi_{\epsilon}^{-1},$ 
$v= (\tau^\prime-t)\varphi_{\epsilon}^{-1}$ and let $\epsilon_1= \min(\epsilon^\prime, \epsilon^{\prime\prime})$, where $\epsilon^\prime= \sup\{\epsilon: \varphi_\epsilon \leq -\frac{c}{A}\}$ and $\epsilon^{ \prime\prime}= \sup\{\epsilon: \varphi_\epsilon \leq -\frac{T-d}{B}\}.$ Then, for $0< \epsilon <\epsilon_1,$
\bea
\;\;\;\\ \nonumber
I_3 &= & \frac{ \epsilon}{\varphi_\epsilon} E_\theta \left[ \int_0^T K\left(\frac{\tau-t}{\varphi_\epsilon} \right) d G_\tau
\right] \\ \nonumber
&=& \frac{ \epsilon}{\varphi_\epsilon} [\int_0^T\int_0^TK\left(\frac{\tau-t}{\varphi_\epsilon} \right)K\left(\frac{\tau^\prime-t}{\varphi_\epsilon}\right)R(\tau,\tau^\prime)d\tau d\tau^\prime]^{1/2}\\\nonumber
&=& \frac{ \epsilon}{\varphi_\epsilon}\varphi_\epsilon[\int_{-\ity}^{\ity} \int_{-\ity}^{\ity}|K(u)K(v) R(t+u \varphi_\epsilon, t+v\varphi_\epsilon)dudv]^{1/2}\\\nonumber
& \leq & C_1\epsilon  (\;\; \mbox{(by using $(A_2)$  )} \\ \nonumber
\eea
for some positive constant $C_1$ by observing that
$$\int_{-\ity}^{\ity} \int_{-\ity}^{\ity}|K(u)K(v) R(t+u \varphi_\epsilon, t+v\varphi_\epsilon)dudv$$
tends to $R(t,t) \int_{-\ity}^{\ity} \int_{-\ity}^{\ity}K(u)K(v)dudv= R(t,t)$
by the condition $(A_1)$ as $\epsilon \raro 0$ by Bochner' theorem (cf. Prakasa Rao (1983)) as $\epsilon \raro 0.$ Since $\epsilon \raro 0$, it follows that $I_1$ tends to zero.
\vsp
Furthermore
\bea
\;\;\;\\ \nonumber
I_2 &= & E_\theta \left[ |\frac{1}{\varphi_\epsilon} \int^T_0 K\left(
\frac{\tau-t}{\varphi_\epsilon}\right) \theta(\tau) x_\tau) d \tau - \theta(t) x_t|
\right] \\ \nonumber
& \leq &  E_\theta \left[ \int^\infty_{-\infty} |K(u)
\left(\theta( t+\varphi_\epsilon u) x_{t+\varphi_\epsilon u}-\theta(t) x_t \right)  \ du| \right]
\\ \nonumber
& \leq &  L\left[ \int^\infty_{-\infty} |K(u) u| \varphi_\epsilon du
\right]\\ \nonumber
& \leq  & C_2 \varphi_\epsilon \\ \nonumber
\eea
for some positive constant $C_2.$  Hence $I_2$ tends to zero as $\epsilon \raro 0.$ Furthermore note that
\bea
\;\;\;\\\nonumber
I_1 &= & E_\theta \left[ |\frac{1}{\varphi_\epsilon} \int^T_0 K
\left(\frac{\tau-t}{\varphi_\epsilon} \right) (\theta (\tau)X_\tau -\theta(\tau) x_\tau)
d\tau |\right] \\\nonumber
&= & E_\theta  \left[ |\int^\infty_{-\infty} K(u)
\left( \theta(t+\varphi_\epsilon u) X_{t+\varphi_\epsilon u}  - \theta(t+\varphi_\epsilon u) x_{t+\varphi_\epsilon u}
)\right) du|\right]\\\nonumber
& \leq &  L E_\theta (\int^\infty_{-\infty} |K(u)|  \left(
|X_{t+\varphi_\epsilon u}-x_{t+\varphi_\epsilon u} |\right) \ du)
\;\;\mbox{(by using the condition $(A_1)$)}\\\nonumber
& \leq &  L E_\theta[\int^\infty_{-\infty} |K(u)|  \sup_{0 \leq t +
\varphi_\epsilon u \leq T}\left(|X_{t+\varphi_\epsilon u}
-x_{t+\varphi_\epsilon u}|\right) \ du ]\\\nonumber
&\leq & L \int^\infty_{-\infty} |K(u)|du \;e^{LT}\epsilon \;E[\sup_{0\leq s \leq T}|G_s|]\\\nonumber
& \leq & C_3 \epsilon \;\;\mbox{(by using (2.4))}\\\nonumber
\eea
for some positive constant $C_4$ depending on $T$ and $L$. Hence $I_3$  tends to zero as  $\epsilon \raro 0.$  Theorem 3.1 is now proved by using the equations (4.1) to (4.4).

\vsp
\NI{\bf Remarks:} From the proof presented above, it is possible to choose the functions $c_\epsilon $ and 
$d_\epsilon $ such that $c_\epsilon \raro 0, d_\epsilon \raro T$ 
and satisfy the conditions 
$$\frac{c_\epsilon}{\varphi_\epsilon}\geq -A, \frac{T-d_\epsilon}{\varphi_\epsilon}\geq B$$
(for instance, choose $c_\epsilon = -A \varphi_\epsilon  $ and $d_\epsilon = T-B \varphi_\epsilon $). Then the estimator $\hat \theta_t$ satisfies the property that
\be
\lim_{\epsilon \rightarrow 0} \sup_{\theta(t) \in \Theta_0(L)} \sup_{c_\epsilon \leq t \leq d_\epsilon } E_\theta ( |\widehat{\theta}_tX_t - \theta(t)x_t|)= 0.
\ee
\vsp
\NI {\bf Proof of Theorem 3.2 :}
Let $ J(t) = \theta(t) x_t.$ By the Taylor's formula, for any $u \in R,$

$$ J(y) = J(u) +\sum^k_{r=1} J^{(r)}(u) \frac{(y-u)^r}{r!} +[ J^{(k)}(z)-J^{(k)}(u)] \frac{(y-u)^k}{k!} $$
for some $z$ such that $|z-u|\leq |y-u|.$ Using this expansion, the equation (4.2) and the conditions  in the expression $I_2$ defined in the proof of  Theorem 3.1, it follows that
\ben
\;\;\\\nonumber
I_2 & \leq &  \left[|\int^\infty_{-\infty} K(u) \left(J(t+\varphi_\epsilon u) - J(t) \right)  \ du |\right]\\\nonumber
&= & [ |\sum^k_{j=1} J^{(j)}(t) (\int^\infty_{-\infty}K(u) u^j du )\varphi^j_\epsilon (j!)^{-1}\\\nonumber
& & \;\;\;\;+(\int^\infty_{-\infty}K(u) u^k (J^{(k)}(z_u) -J^{(k)}(x_t))du \varphi^k_\epsilon (k !)^{-1}|]\\ \nonumber
\een
for some $z_u$ such that $|x_t-z_u|\leq |x_{t+\varphi_\epsilon u}-x_t| \leq C|\varphi_\epsilon u|$ for some positive constant $C.$ Hence
\bea
I_2 & \leq & C_4 L \left[  \int^\infty_{-\infty} |K(u)u^{k+1}|\varphi^{k+1}_\epsilon (k!) ^{-1}  du  \right]\\\nonumber
& \leq & C_5 (k!)^{-1} \varphi^{k+1}_\epsilon \int^\infty_{-\infty} |K(u) u^{k+1}| du  \\\nonumber
&\leq & C_6 \varphi_\epsilon^{k+1}\\ \nonumber
\eea
for some positive constant $C_6$ depending on $A, B,T$ and $L$. Combining the relations (4.2), (4.4) and (4.6) , we get that there exists a positive constant $C_7$ depending on $T,L,A, B$ such that 
$$ \sup_{c \leq t \leq d}E_\theta|\widehat{\theta}_tx_t-\theta(t) x_t| \leq C_{7} (\epsilon+  \varphi^{k+1}_\epsilon +\epsilon). $$ 
Choosing $ \varphi_\epsilon = \epsilon^{\frac{1}{k+1}},$  we get that 
$$ \limsup_{\epsilon \rightarrow 0} \sup_{\theta(.) \in \Theta_{k+1}(L) } \sup_{c \leq t \leq d} E_\theta |\widetilde{\theta}_t X_t - \theta(t) x_t|\epsilon^ {-1} < \infty. $$ 
This completes the proof of Theorem 3.2. 
\vsp 
\vsp \NI{\bf Proof of Theorem 3.3:}

From the equation (3.1), we obtain that
\bea
\lefteqn{ \widehat{\theta}_t X_t -\theta(t) x_t}\\\nonumber
 &= &[ \frac{1}{\varphi_\epsilon} \int^T_0 K \left(\frac{\tau-t}{\varphi_\epsilon} \right)
 \left( \theta(\tau) X_\tau- \theta(\tau) x_\tau \right) \  d \tau \\\nonumber
 & & + \frac{1}{\varphi_\epsilon} \int^T_0 K \left( \frac{\tau-t}{\varphi_\epsilon}\right) \theta(\tau) x_\tau d\tau -\theta(t)x_t+ \frac{\epsilon}{\varphi_\epsilon} \int^T_0 K \left( \frac{\tau-t}{\varphi_\epsilon}\right) dG_\tau]\\\nonumber
 &= & [ \int^\infty_{-\infty} K(u) (\theta(t+\varphi_\epsilon u) X_{t+\varphi_\epsilon u} - \theta (t+\varphi_\epsilon u) x_{t+\varphi_\epsilon u}) \ du  \\\nonumber
 & & +\int^\infty_{-\infty} K(u) (\theta (t+\varphi_\epsilon u) x_{t+\varphi_\epsilon u}- \theta(t) x_t) \ du \\\nonumber
 & &+ \frac{\epsilon}{\varphi_{\epsilon}}\int^T_0 K\left(\frac{\tau-t}{\varphi_\epsilon}
 \right) d G_\tau].\\\nonumber
\eea
Let $J(t)= \theta(t)x_t.$ By the Taylor's formula, for any $u \in R,$
$$ J(y) = J(u) +\sum^{k+1}_{r=1} J^{(r)}(u) \frac{(y-u)^j}{j!} +[J^{(k+1)}(z)-J^{(k+1)}(x)] \frac{(y-u)^{k+1}}{(k+1)!} $$
for some $z$ such that $|z-u|\leq |y-u|.$ Let 
$$m= \frac{J^{(k+1)}(x_t)}{(k+1)!}\int_{-\ity}^{\ity}K(u)u^{k+1}du$$
and
$$ R_1(t)= \varphi^{-(k+1)}_\epsilon\int^T_0 K\left(\frac{\tau-t}{\varphi_\epsilon} \right) (\theta(\tau) X_\tau -\theta (\tau) x_\tau)
d\tau.$$
By arguments similar to those given in (4.3) for obtaining upper bounds, it follows that 
$$E|R_1(t)| \leq C \varphi_\epsilon^{-k}\epsilon.$$
Let 
$$ R_2(t)= \varphi_\epsilon^{-(k+1)} \int^T_0 K\left(\frac{\tau-t}{\varphi_\epsilon} \right) (\theta(\tau) x_\tau -\theta (t) x_t)
d\tau.$$
Observe that
$$R_2(t)=m+o(1)$$
by an application of the Taylor's expansion under the condition $(A_3).$ \\
\vsp
Furthermore
\bea
\varphi_\epsilon^{-(k+1)}(\hat \theta_t x_t-\theta(t) x_t) &= & \epsilon \varphi^{-(k+2)}_\epsilon\int^T_0 K\left(\frac{\tau-t}{\varphi_\epsilon} \right) d G_\tau + R_2(t) + R_1(t)\\\nonumber
&=&  \epsilon \varphi^{-(k+2)}_\epsilon\int^T_0 K\left(\frac{\tau-t}{\varphi_\epsilon} \right) d Z_\tau + m+o(1)+O_p(\varphi_\epsilon^{-k}\epsilon).\\\nonumber 
\eea
Let $\varphi_\epsilon$ be chosen so that $(\varphi_\epsilon)^{-1}= \epsilon \varphi_\epsilon^{-(k+2)}.$ Such a choice is $\varphi_\epsilon =\epsilon^v$ where $v= (k+1)^{-1}.$ 
We will now study the asymptotic behaviour of the random variable 
$$W_\epsilon= (\varphi_\epsilon)^{-1}\int^T_0 K\left(\frac{\tau-t}{\varphi_\epsilon} \right) d Z_\tau $$
as $\epsilon \raro 0.$  Note that	
\be
\varphi_\epsilon^{-(k+1)}(\hat \theta_t X_t-\theta(t) x_t)=W_\epsilon+ m+o_p(1)
\ee
Note that $W_\epsilon$ is Gaussian with mean zero and variance $R(t,t)+o(1)$ as $\epsilon \raro 0.$
$$m= \frac{J^{(k+1)}(x_t)}{(k+1)!}\int_{-\ity}^{\ity}K(u)u^{k+1}du$$
where $J(t)=\theta(t)x_t.$
Hence 
$$\varphi_\epsilon^{-(k+1)}(\hat \theta_t X_t-\theta(t) x_t)$$
is Gaussian with mean 
$$m= \frac{J^{(k+1)}(x_t)}{(k+1)!}\int_{-\ity}^{\ity}K(u)u^{k+1}du$$
and variance $R(t,t)$ as $\epsilon \raro 0.$ 
\vsp
Note that the results given above deal with asymptotic  properties of the estimator for the function 
$$J(t)= \theta(t)x_t= \theta(t) x_0 \exp(\int_0^t \theta(s)\;ds).$$ 
\vsp
We will now present another method for the estimation of the linear multiplier $\theta(t).$
\vsp
\newsection{Estimation of the Multiplier $\theta(.)$}

Let $\Theta_\rho(L_\gamma)$ be a class of functions uniformly bounded  and $k$-times continuously differentiable for some integer 
$k \geq 1 $ with the $k$-th derivative satisfying the Holder condition of the order $\gamma \in (0,1):$
$$|\theta^{(k)}(t)-\theta^{(k)}(s)|\leq L_\gamma |t-s|^\gamma, \rho=k+\gamma.$$

From the Lemma 3.1, it follows that
$$|X_t-x_t| \leq \epsilon e^{Lt}\sup_{0\leq s \leq T} |G_s|.$$
Let
$$A_t= \{\omega: \inf_{0\leq s \leq t}X_s(\omega)\geq \frac{1}{2}x_0e^{-Lt}\}$$
and let $A=A_T.$ Define the process $Y$ with the differential
$$dY_t=\theta(t) I(A_t) dt + \epsilon X_t^{-1}I(A_t)\;dG_t, 0\leq t \leq T.$$
We will now construct an estimator of the function $\theta(.)$ based on the observation of the process $Y$ over the interval $[0,T].$
Define the estimator
$$\tilde \theta(t)= I(A) \frac{1}{\varphi_\epsilon}\int_0^T K(\frac{t-s}{\varphi_\epsilon})dY_s$$
where the kernel function $K(.)$ satisfies the conditions $(A1)-(A3)$. Observe that
\ben
E|\tilde \theta(t)-\theta(t)| &= & E|I(A) \frac{1}{\varphi_\epsilon}\int_0^T K(\frac{t-s}{\varphi_\epsilon})(\theta(s)-\theta(t))ds\\\nonumber
&& \;\;\;\; + I(A^c)\theta(t)+I(A)\frac{\epsilon}{\varphi_\epsilon}\int_0^T K(\frac{t-s}{\varphi_\epsilon})X_s^{-1}dG_s|\\\nonumber
&\leq & E|I(A)\int_R K(u)[\theta(t+u\varphi_\epsilon)-\theta(t)]du|+ |\theta(t)| P(A^c)\\\nonumber
&& \;\;\;\; + \frac{\epsilon}{\varphi_\epsilon}|E[I(A)\int_0^T K(\frac{t-s}{\varphi_\epsilon})X_s^{-1}dG_s]|\\\nonumber
&=& I_1+I_2+I_3. \;\;\mbox{(say)}.\\\nonumber
\een
Applying the Taylor's theorem and using the fact that the function $\theta(t)\in \Theta_\rho(L_\gamma)$, it follows that
\ben
I_1 \leq \frac{L_\gamma}{(k+1)!}\varphi_\epsilon^\rho \int_R|K(u)u^\rho|du.
\een
Note that, by Lemma 3.1,
\ben
P(A^c) & = & P(\inf_{0\leq t \leq T}X_t < \frac{1}{2}x_0e^{-LT})\\\nonumber
&\leq & P(\inf_{0 \leq t \leq T}|X_t-x_t| + \inf_{0\leq t \leq T}x_t < \frac{1}{2}x_0e^{-LT})\\\nonumber
&\leq & P(\inf_{0 \leq t \leq T}|X_t-x_t| < -\frac{1}{2}x_0e^{-LT})\\\nonumber
&\leq & P(\sup_{0 \leq t \leq T}|X_t-x_t| > \frac{1}{2}x_0e^{-LT})\\\nonumber
&\leq & P(\epsilon e^{LT}\sup_{0 \leq t \leq T}|G_t|>\frac{1}{2}x_0e^{-LT})\\\nonumber
&= & P(\sup_{0 \leq t \leq T}|G_t|>\frac{x_0}{2\epsilon}e^{-2LT})\\\nonumber
&\leq &(\frac{x_0}{2\epsilon}e^{-2LT})^{-1} E{\sup_{0 \leq t \leq T}|G_t|}\\\nonumber
&\leq & D\epsilon
\een
for some positive constant $D$ by the assumption on the Gaussian process $G$. The upper bound obtained above and the fact that $|\theta(s)|\leq L, 0\leq s \leq T$ leads  an upper bound for the term $I_2.$ We have used the inequality
$$x_t= x_0 \exp(\int_0^t\theta (s)ds)\geq x_0 e^{-Lt}$$
in the computations given above. Furthermore
\ben
\lefteqn{|E[I(A)\int_0^T K(\frac{t-s}{\varphi_\epsilon})X_s^{-1}dG_s|}\\\nonumber
&=& |E[\int_0^T K(\frac{t-s}{\varphi_\epsilon})X_s^{-1}I(A_s)dG_s|\\\nonumber
&\leq & C_9 E[(\int_0^T|K(\frac{t-s}{\varphi_\epsilon})X_s^{-1}I(A_s)| dG_s)\\\nonumber
&\leq &C_{10} e^{LT} (\int_0^T|G(\frac{t-s}{\varphi_\epsilon})|dG_s)\\\nonumber
& \leq & C_{11} \varphi_\epsilon (R(t,t)+o(1))\\\nonumber
\een
for some positive constant $C_{11}$ which leads to an upper bound on the term $I_3.$ Combining the above estimates, it follows that
\ben
E|\tilde \theta(t)-\theta(t)|\leq C(\varphi_\epsilon^\rho +  \epsilon + \epsilon) 
\een
for some positive constant $C.$ Choosing $\varphi_\epsilon=\epsilon^{1/\rho},$ we obtain that
\ben
E|\tilde \theta(t)-\theta(t)|\leq C\epsilon
\een
for some positive constants $C$. Hence we obtain the following result implying the uniform consistency of the estimator $\tilde \theta(t)$ as an estimator of $\theta(t)$ as $\epsilon \raro 0.$\\

\NI{\bf Theorem 5.1:} {\it Let $\theta \in \Theta_\rho(L)$ and $\varphi_\epsilon= \epsilon^{1/\rho}.$ Suppose the conditions $(A_1)-(A_3)$ hold. Then, for any interval $[c,d] \subset [0,T],$ }
\ben
\limsup_{\epsilon \raro 0}\sup_{\theta(.)\in \Theta_\rho(L)}\sup_{c\leq t \leq d}E|\tilde \theta(t)-\theta(t)| \epsilon^{-1}<\ity.
\een
\vsp
\NI{\bf Funding :} This work was supported under the scheme ``INSA Senior Scientist" by the Indian National Science Academy  while the author was at the CR Rao Advanced Institute for Mathematics, Statistics and Computer Science, Hyderabad 500046, India.\\
\vsp
\NI{\bf References :}
\begin{description}

\item Chen, Y., and Zhou, H. (2020) Parameter  estimation for an Ornstein-Uhlenbeck process driven by a general Gaussian noise, arxiv:2002.09641v1 [math.PR] 22 Feb 2020.

\item El Machkouri, M., Es-Sebaiy, K., and Ouknine, Y. (2015) Parameter estimation for the non-ergodic Ornstein-Uhlenbeck processes driven by Gaussian process, arXiv:1507.00802v1 [math.PR] 3 July 2015.

\item Kutoyants, Yu. (1994) {\it Identification of dynamical Systems with Small Noise}, Kluwer, Dordrecht.

\item Lu, Y. (2022) Parameter estimation of non-ergodic Ornstein-Uhlenbeck processes driven by general Gaussian processes, arXiv:2207:13355v1 [math.ST] 27 Jul 2022. 

\item Mishura, Y. and Zili, M. (2018) {\it Stochastic Analysis of Mixed Fractional Gaussian Processes}, ISTE Press and Elsevier, UK.

\item Nourdin, I. (2012) {\it Selected Aspects of Fractional Brownian Motion}, Bocconi and Springer Series, Bocconi University Press, Milan.

\item Prakasa Rao, B.L.S. (1983) {\it Nonparametric Functional Estimation}, Academic Press, New York.

\item  Prakasa Rao, B.L.S. (1999)  {\it Statistical Inference for Diffusion Type Processes}, Arnold, London and Oxford University Press, New
York.
 
\item Prakasa Rao, B.L.S. (2010) {\it Statistical Inference for Fractional Diffusion Processes}, Wiley, London.

\item  Prakasa Rao, B.L.S. (2021) Nonparametric estimation of linear multiplier in stochastic differential equations driven by $\alpha$-stable noise, {\it Journal of the Indian Statistical Association}, {\bf 59}, 65-81.

\item Prakasa Rao, B.L.S. (2022) Minimum $L_1$l-norm estimation for fractional ornstein-Uhlenbeck process driven by a gaussian process, arxiv:2208.04366v1 [math. PR] 8 Aug 2022.
\end{description}
\vsp
\NI {CR Rao Advanced Institute of Mathematics, Statistics and Computer Science, Hyderabad, India.}\\ 
\NI{e-mail: blsprao@gmail.com}
\end{document}